\documentclass[preprint,review,12pt]{amsart}
\usepackage{amsmath,amscd,amssymb,amsfonts,amsthm}
\usepackage[all]{xy}
\xyoption{all}

\newtheorem*{thm*}{Theorem}
\newtheorem*{cor*}{Corollary}
\newtheorem{thm}{Theorem}[section]

\newtheorem{cor}[thm]{Corollary}
\newtheorem{lem}[thm]{Lemma}

\theoremstyle{definition}
\newtheorem{defn}[thm]{Definition}

\newcommand{\C}{{\mathbb{C}}}

\newcommand{\SH}[1][]{{\mathcal{S}_{#1}H}}

\newcommand{\SQ}{{\mathcal{S}Q}}
\newcommand{\SG}{{\mathcal{S}G}}
\newcommand{\SR}{{\mathcal{S}R}}

\newcommand{\SGH}{{\mathcal{S}G/H}}

\newcommand{\N}{{\mathbb{N}}}
\newcommand{\I}{{\mathcal{I}}}
\newcommand{\tensor}{{\otimes}}

\newcommand{\btensor}{{\hat\otimes}}

\newcommand{\isom}{{\, \cong \,}}

\newcommand{\Bdd}{\mathfrak{B}}

\newcommand{\im}{\operatorname{im}}
\newcommand{\Hom}{\operatorname{Hom}}
\newcommand{\Ext}{\operatorname{Ext}}
\newcommand{\bHom}{\operatorname{bHom}}
\newcommand{\bExt}{\operatorname{bExt}}

\begin{document}


\title{A Spectral Sequence for Polynomially Bounded Cohomology}

\author{Bobby Ramsey}
\address{The Ohio State University\\Dept. of Mathematics\\321 W. 18th Ave.\\Columbus, OH  43210}

\begin{abstract}
We construct an analogue of the Lyndon-Hochschild-Serre spectral sequence in the context of
polynomially bounded cohomology.  For $G$ an extension of $Q$ by $H$, this spectral sequences converges to 
the polynomially bounded cohomology of $G$, $HP^*(G)$.  If the extension is a polynomial 
extension in the sense of Noskov with $H$ and $Q$ isocohomological and $Q$ of type $HF^\infty$, the spectral sequence 
has $E_2^{p,q}$-term $HP^q(Q; HP^p(H))$, and $G$ is isocohomological for $\C$.  By referencing results of 
Connes-Moscovici and Noskov if $H$ and $Q$ are both isocohomological and have the Rapid Decay property, then
$G$ satisfies the Novikov conjecture.
\end{abstract}


\maketitle
\section{Introduction}

In \cite{CM} Connes and Moscovici prove the Novikov conjecture for all finitely generated discrete 
groups satisfying two properties.
The first is the Rapid Decay property of Jolissaint \cite{Jol}, which ensures the
existence of a smooth dense subalgebra of the reduced group $C^*$-algebra.
The second property is that every cohomology class can be represented by a cocycle of 
polynomial growth, with respect to some (hence any) word-length function on the group.  

The polynomially bounded cohomology of a group $G$, denoted $HP^*(G)$, obtained
by considering only cochains of polynomial growth, has been of interest recently.  
The inclusion of these polynomially bounded cochains into the full cochain complex yields a homomorphism
from the polynomially bounded cohomology to the full cohomology of the group.  The
second property of Connes-Moscovici above is that this polynomial comparison
homomorphism is surjective.  A group $G$ is {\it{isocohomological for $M$}} if
$HP^*(G; M)$ is bornologically isomorphic to $H^*(G;M)$.  The term \emph{isocohomological} is taken from Meyer 
\cite{M4}, where it describes a homomorphism between two bornological
algebras.  What is meant here by `$G$ is isocohomological for $\C$', is a weakened version of 
what Meyer refers to as \emph{the embedding $\C[G] \to \SG$ is isocohomological},
where $\SG$ refers to the Fr\'echet algebra of functions $G \to \C$ of $\ell^1$-rapid decay.  
We call $G$ \emph{isocohomological} if it is isocohomological for all $SG$ modules.

Influenced by Connes and Moscovici's approach, in \cite{Ji1} Ji defined polynomially bounded cohomology and showed 
that virtually nilpotent groups,
are isocohomological for $\C$.  In \cite{M1} Meyer showed that polynomially combable groups are isocohomological for $\C$.
By citing a result of Gersten regarding classifying spaces for
combable groups \cite{Ger}, Ogle independently showed that polynomially combable groups are isocohomological for $\C$.
In \cite{JR} it was shown that for the class of finitely presented FP$^\infty$ groups, a group is isocohomological
if and only if it has Dehn functions which are polynomially bounded in all dimensions.

For a group extension, $0 \to H \to G \to Q \to 0$,
the Lyndon-Hochschild-Serre (LHS) spectral sequence is a first-quadrant spectral sequence with
$E_2$-term $H^*( Q; H^*( H ) )$ and converging to $H^*(G)$.
In \cite{N2} Noskov generalized the construction of the LHS spectral sequence to obtain a 
spectral sequence in bounded cohomology for which, under suitable topological circumstances, one could
identify the $E_2$-term as $H_b^*( Q; H_b^*( H ) )$ and which converges to $H_b^*( G )$.
Ogle considered the LHS spectral sequence in the context of P-bounded cohomology in \cite{O}.
There additional technical considerations are also needed to ensure the appropriate $E_2$-term.  It is
not clear for which class of extensions these conditions are satisfied.

In this paper we resolve this issue for polynomial extensions, which were proposed by Noskov 
in \cite{N}.  
Let $\ell$, $\ell_H$, and $\ell_Q$ be word-length functions on $G$, $H$, and $Q$ respectively, and let
\[ 0 \to H \to G \stackrel{\pi}{\to} Q \to 0 \]
be an extension of $Q$ by $H$.  
Let $q \mapsto \overline q$ be a cross section of $\pi$.
To this cross section there is associated a function
$[ \cdot, \cdot ] : Q \times Q \to H$ by
$\overline q_1 \overline q_2 = \overline{q_1 q_2} [ q_1, q_2 ]$, called the \emph{factor set} 
of the extension.  The factor set has polynomial
growth if there exists constants $C$ and $r$ such that 
$\ell_H ([q_1, q_2])  \leq C( (1 + \ell_Q(q_1))( 1 + \ell_Q(q_2)) )^r$.
The cross section also determines a \emph{set-theoretical action} of $Q$ on $H$.  This
is a map $Q \times H \to H$ given by
$(q,h) \mapsto h^q = \overline{q}^{-1} h \overline{q}$.  
The set-theoretical action is polynomially bounded if
there exists constants $C$ and $r$ such that
$\ell_H( h^q ) \leq C \ell_H(h)  ( 1 + \ell_Q(q) )^r$.
In what follows, we adopt the convention that if $\mathcal{Q}$ is the finite generating set for 
$Q$ and $\mathcal{A}$ is the finite
generating set for $H$, then as the generating set for $G$ we will take
the set of $h \in \mathcal{A}$ and $\overline q$ for $q \in \mathcal{Q}$.

\begin{defn}
An extension $G$ of a finitely generated group $Q$ by a finitely generated group $H$ is said to 
be a \underline{polynomial extension} if there is a cross section yielding a factor set of 
polynomial growth and inducing a polynomial set-theoretical action of $Q$ on $H$.
\end{defn}
Our main theorem is as follows.

\begin{thm*}
Let $0 \to H \to G \rightarrow Q \to 0$ be a polynomial extension of the $HF^\infty$ group $Q$,
with both $H$ and $Q$ isocohomological.
There is a bornological spectral sequence with $E^{p,q}_2 \isom HP^p( Q; HP^q(H) )$ 
which converges to $HP^*(G)$.
\end{thm*}

We can compare this spectral sequence with the LHS spectral sequence.

\begin{cor}\label{isoExtension}
Let $0 \to H \to G \to Q \to 0$ be a polynomial extension with $Q$ of type $HF^\infty$.
If $H$ and $Q$ are isocohomological, then $G$ is isocohomological for $\C$.
\end{cor}

Applying a result of Noskov \cite{N} regarding the polynomial
extension of groups with the Rapid Decay property, as well as the results of Connes-Moscovici,
we obtain the following corollary.
\begin{cor}\label{novExtension}
Let $0 \to H \to G \to Q \to 0$ be a polynomial group extension with $H$ and $Q$ isocohomological, 
and $Q$ of type $HF^\infty$.
If both $Q$ and $H$ have the Rapid Decay property, $G$ satisfies the Novikov conjecture.
\end{cor}

It would be convenient to work in the category of Fr\'echet and DF spaces, 
however there are many quotients involved in the construction yielding spaces which 
need be neither Fr\'echet nor DF.  This is the issue Ogle overcomes by use of a technical
hypothesis in \cite{O}, and it is at the heart of the topological consideration in \cite{N2}, used to 
identify the $E_2$-terms in their spectral sequences.
To overcome these obstacles we work mostly in the bornological category
and utilize an adjointness relationship of the form 
\[ \Hom( A \btensor B, C ) \isom \Hom( A, \Hom( B, C ) )\] which is not true in the
category of locally convex topological vector spaces.  A bornology on a space
is an analogue of a topology, in which boundedness replaces openness as the key consideration.
In this context, we are also able to bypass many of the issues involved in the topological analysis of
vector spaces. When endowed with the fine bornology, as defined later, any complex vector space is a
complete bornological vector space.  The finest topology yielding a complete topological structure
on such a space is cumbersome.  This bornology allows us to replace analysis of continuity in this
topology, to boundedness in finite dimensional vector spaces.

In Section \ref{Born}, we recall the relevant concepts from the bornological framework developed by 
Hogbe-Nlend and extended by Meyer and others.  Section \ref{bSS} consists of translating the usual algebraic 
spectral sequence arguments into this framework.  This is mainly verifying that the vector space isomorphisms 
are in fact isomorphisms of bornological spaces.  These will be the main tools of our construction.
In Section \ref{pcoh} we define the relevant bornological algebras and
define the polynomially bounded cohomology of a discrete group endowed with a length function, as
well as basic materials for the construction of our spectral sequence and some of its applications.
The actual computation is the focus of the final section. 

I wish to thank Ron Ji, Crichton Ogle, Ralf Meyer and the Referee for their invaluable comments.  
This article represents a revised version of my thesis.

\section{Bornologies}\label{Born}

Let $A$ and $B$ be subsets of a locally convex topological vector space $V$.  A subset $A$ is {\it circled} if 
$\lambda A \subset A$ for all $\lambda \in \C$, $|\lambda| \leq 1$.  It is a {\it disk} if
it is both circled and convex.  For two subsets, $A$ {\it absorbs} $B$ if there is an $\alpha > 0$ such that
$B \subset \lambda A$ for all $| \lambda | \geq \alpha$, and  $A$ is {\it absorbent} if it absorbs
every singleton in $V$.  The {\it circled (disked, convex) hull} of $A$ is the smallest circled 
(disked, convex) subset of $V$ containing $A$.
For an absorbent set $A$ there is associated a semi-norm on $\rho_A$ on $V$
given by $\rho_A( v ) = \inf \left\{ \alpha > 0 ; v \in \alpha A \right\}$.  For an arbitrary
subset $A$, denote by $V_A$ the subspace of $V$ spanned by $A$.  If $A$ is a disked set,
$\rho_A$ is a semi-norm on $V_A$.  $A$ is a {\it completant disk} if $V_A$ is a Banach space in the 
topology induced by $\rho_A$.

Let $V$ be a locally convex topological vector space.  A convex vector space bornology on $V$ is a collection
$\Bdd$ of subsets of $V$ such that:
\begin{enumerate}
\item For every $v \in V$, $\{v\} \in \Bdd$.
\item If $A \subset B$ and $B \in \Bdd$ then $A \in \Bdd$.
\item If $A$, $B \in \Bdd$ and $\lambda \in \C$, then $A + \lambda B \in \Bdd$.
\item If $A \in \Bdd$, and $B$ is the disked hull of $A$, then $B \in \Bdd$.
\end{enumerate}
Elements of $\Bdd$ are said to be the bounded subsets of $V$.

An important example in what follows is the {\it{fine bornology}}.  A set is bounded in
the fine bornology on $V$ if it is contained and bounded in some finite dimensional subspace of $V$.

Let $V$ and $W$ be two bornological spaces.  A map $V \to W$ is {\it bounded} if 
the image of every bounded set in $V$ is bounded in $W$.  A {\it bornological isomorphism} is a
bounded bijection with bounded inverse.  The collection of all bounded linear maps $V \to W$
is denoted $\bHom(V,W)$.  There is a canonical complete convex bornology on $\bHom(V,W)$, given
by the families of equibounded functions;  a family $\mathcal{U}$ is equibounded if
for every bounded $A \subset V$, $\mathcal{U}(A)=\{ u(a) \, | \, u \in U, a \in A\}$ is bounded in $W$.

If $W$ is a bornological space and $V \subset W$, then there is a bornology on $V$ induced
from that on $W$ in the obvious way.  There is also a bornology on $W/V$ induced from $W$.
A subset $B \subset W/V$ is bounded if and only if there is a bounded $C \subset W$
which maps to $B$ under the canonical projection $W \to W/V$.

Let $V$ be a bornological vector space.  A sequence $(v_i)$ in $V$ {\it converges bornologically}
to $0$ if there is a bounded subset $B \subset V$ and a sequence of scalars $\lambda_i$ tending to 
$0$ such that for all $i$, $v_i \in \lambda_i B$.  $(v_i)$ converges bornologically to $v$ if
$(v_i - v)$ converges bornologically to $0$.

Let $V$ and $W$ be bornological vector spaces.  The complete projective bornological tensor product
$V \btensor W$ is given by the following universal property:
For any complete bornological vector space $X$, a jointly bounded bilinear map
$V \times W \to X$ extends uniquely to a bounded map $V \btensor W \to X$.
Unlike the complete topological projective tensor product on the category of
locally compact vector spaces, this bornological tensor product admits an adjoint.

\begin{lem}[\cite{M4}]
Let $A$ be a complete bornological algebra, $V$ and $M$ complete bornological $A$-modules,
and $W$ a complete bornological space.  There is a bornological isomorphism
$\bHom_{A}( V \btensor W, M ) \isom \bHom( W, \bHom_{A}( V,M ) )$.
\end{lem}

We will be interested in bornologies on Fr\'echet spaces.  For a Fr\'echet space $F$,
there is a countable directed family of seminorms, $\| \cdot \|_n$ yielding the topology.
A set $U$ is said to be von Neumann bounded if $\| U \|_n < \infty$ for all $n$.
The collection of all von Neumann bounded sets forms the von Neumann Bornology on $F$.
This will be our bornology of choice on Fr\'echet spaces, due to the relation with
topological constructs.

\begin{defn}\label{defnNet}
A \underline{net} in a bornological vector space $F$ is a family of disks, $\{e_{i_1, i_2, \ldots, i_k}\}$ in $F$, 
indexed by $\bigcup_{k \in \N} \I^k$ ( $\I$ some countable set ), which satisfies the following conditions.
\begin{enumerate}
\item $F = \bigcup_{i \in \I} e_i$, and 
	$e_{i_1, \ldots, i_k} = \bigcup_{i \in \I} e_{i_1, \ldots, i_k, i}$ for $k > 1$.
\item For every sequence $(i_k)$ in $\I$, there is a sequence $(\nu_k)$ of positive reals such that
	for each $f_k \in e_{i_1, \ldots, i_k}$ and each $\mu_k \in [ 0, \nu_k ]$, the series
	$\sum_{k = 1}^{\infty} \mu_k f_k$ converges bornologically in $F$, and for each $k_0 \in \N$, the
	series $\sum_{k = k_0}^{\infty} \mu_k f_k$ lies in $e_{i_1, \ldots, i_{k_0}}$.
\item For every sequence $(i_k)$ in $\I$ and every sequence $(\lambda_k)$ of positive reals,
	$\bigcup_{k=1}^{\infty} \lambda_k e_{i_1, \ldots, i_k}$ is bounded in $F$.
\end{enumerate}
\end{defn}

As an example, Hogbe-Nlend shows that every bornological space with a countable
base has a net \cite[p58]{H}.

\begin{lem}\label{subspaceNet}
Let $F$ be a bornological vector space with a net, and let $V$ be a subspace of $F$.  
Then $V$ has a net.
\end{lem}

\begin{lem}\label{quotientNet}
Let $F$ be a bornological vector space with a net, and let $V$ be a subspace of $F$.
Then $F/V$ has a net.
\end{lem}
\begin{proof}
Denote the net on $F$ by $\mathcal{R} = \left\{e_{i_1, \ldots, i_k}\, | \, i_j \in \I\right\}$,  and
let $\pi: F \to F/V$ be the projection.  
Set $\mathcal{R}' = \left\{e'_{i_1, \ldots, i_k} = \pi e_{i_1, \ldots, i_k}\right\}$.  
\end{proof}

\begin{lem}\label{homNet}
Let $U$ be a Fr\`echet space.  Then $\bHom(U,\C^N)$ has a net.
\end{lem}
\begin{proof}
In this case, boundedness and continuity are equivalent for homomorphisms,
and the equibounded families are precisely the equicontinuous families.

For each finite sequence of ordered triples of positive integers $(n_1,M_1,K_1)$, $\ldots$, $(n_k, M_k, K_k)$, 
define $b_{ (n_1,M_1, K_1), \ldots, (n_k,M_k, K_k)}$ to be the set of all $f \in \bHom(U,\C^N)$ such that for all 
$i$ between $1$ and $k$, $|f(u)| < M_i$ for all $u \in U$ with $\| u \|_{U,n_i} < K_i$.  We show that
this gives a countable base for the bornology on $\bHom(U,\C^N)$.  

If $W$ is an equibounded family, then for all neighborhoods $V$ of zero in $\C^N$,
there exist $n_1, \ldots, n_k$ and $K_1, \ldots, K_k$  such that 
$\{ u \in U \, | \, \text{ for all } 1 \leq i \leq k \, \| u \|_{U,n_i} < K_i \}$ is contained
in $W^{-1}(V)$.  For each $f \in W$ and for each $u$ in this set, $f(u) \in V$.
Let $V$ be the open ball of radius $1$ in $\C^N$, and let $n_1, \ldots, n_k$ and $K_1, \ldots, K_k$
be as above.  Then 
$W \subset b_{ (n_1,1, K_1), \ldots, (n_k, 1, K_k)}$.
Let $B=b_{ (n_1,M_1, K_1), \ldots, (n_k,M_k, K_k)}$, $f \in B$, $M = \max M_i$, and let 
$V$ be the open ball of radius $R$ in $\C$.  Then $V = \frac{R}{M} V'$ where $V'$ is the open
ball of radius $M$ in $\C^N$.  $f^{-1}(V) = \frac{R}{M}f^{-1}(V')$, so if $B^{-1}(V')$ is a neighborhood
of zero in $U$, then so is $B^{-1}(V)$.
Let $u \in U$ be such that for all $1 \leq i \leq k$ we have $\| u \|_{U,n_i} < K_i$.
Then $|f(u)| < M$, so $f(u) \in V'$.  Let $S$ be the set of all such $u \in U$.  It
is a neighborhood of zero.  Moreover $f(S) \subset V'$ so that $S \subset f^{-1}(V')$,
whence $S \subset B^{-1}(V')$.
This implies that the bornology on $\bHom(U,\C^N)$ has a countable base.
\end{proof}

Our interest in nets is the following analogue of the open-mapping theorem.
\begin{thm}[{ \cite[p61]{H} }]\label{bddMapIsom} 
Let $E$ and $F$ be convex bornological spaces such that $E$ is complete and $F$ has a
net.  Every bounded linear bijection $v:F \to E$ is a bornological
isomorphism.
\end{thm}

\section{Preliminary Results on Bornological Spectral Sequences}\label{bSS}
This section contains several results from McCleary's book \cite{Mc} translated into the bornological
framework.  The proofs given follow McCleary, with modifications to verify that the vector space
isomorphisms involved are isomorphisms of bornological spaces.

Let $(A,d)$ be a differential graded bornological module.  That is, $A=\bigoplus_{n=0}^{\infty} A^n$ is a graded 
bornological module and $d : A \to A$ is a degree $1$ bounded linear map with $d^2 = 0$.   Let $F$
be a filtration of $A$ which is preserved by the differential, so that for all $p$, $q$ we have 
$d(F^pA^q) \subset F^pA^{q+1}$.  Assume further that the filtration is decreasing, in that
$ \ldots \subset F^{p+1}A^q \subset F^pA^q \subset F^{p-1}A^q \subset \ldots $.
Such an $(A,F,d)$ will be referred to as a \emph{filtered differential graded bornological module}.
Denote by $d^{p,q} : F^pA^{p+q} \to F^pA^{p+q+1}$ the restriction of $d$, so $d$ is the direct
sum of $d^{p,q}$. 
The filtration $F$ is said to be bounded if for each $n$, there is $s = s(n)$ and $t = t(n)$ such that 
\[ 0 = F^sA^n \subset F^{s-1}A^n \subset \ldots \subset F^{t+1}A^n \subset F^tA^n = A^n \]
Let
\begin{eqnarray*}
Z_r^{p,q} & = & F^pA^{p+q} \cap (d^{p+r,q-r})^{-1}(F^{p+r}A^{p+q+1})\\
B_r^{p,q} & = & F^pA^{p+q} \cap d^{p-r,q+r-1}(F^{p-r}A^{p+q-1})\\
Z_\infty^{p,q} & = & F^pA^{p+q} \cap \ker d \\
B_\infty^{p,q} & = & F^pA^{p+q} \cap \im d
\end{eqnarray*}
where each of these subspaces are given the subspace bornology.
Let $d^n : A^n \to A^{n+1}$ be the restriction of $d$.
These definitions yield the following `tower' of submodules.
\[ B_0^{p,q} \subset B_1^{p,q} \subset \ldots \subset B_\infty^{p,q} \subset` Z_\infty^{p,q}
 \subset \ldots \subset Z_1^{p,q} \subset Z_0^{p,q} \]
Moreover $d^{p-r,q+r-1}( Z_r^{p-r,q+r-1} ) = B_r^{p,q}$.

If the filtration is bounded and $r \geq \max\{ s(p+q+1)-p, p-t(p+q-1)\}$ then 
$(d^{p+r,q-r})^{-1}(F^{p+r}A^{p+q+1})$ is the kernel of $d$, 
$Z_r^{p,q} = Z_\infty^{p,q}$, and $B_r^{p,q} = B_\infty^{p,q}$.

\begin{lem}\label{filtDiffGrdBornModSS}
For $(A,F,d)$ a filtered differential graded bornological module, there is a spectral sequence
of bornological modules
$(E^{*,*}_r, d_r)$, $r = 1$, $2$, $\ldots$, with $d_r$ of bidegree $(r, 1-r)$ and
$E_1^{p,q} \isom H^{p+q}(F^{p}A/F^{p+1}A)$.  If the filtration is bounded the
spectral sequence converges to $H(A,d)$, $E_{\infty}^{p,q} \isom F^pH^{p+q}(A,d)/F^{p+1}H^{p+q}(A,d)$.
\end{lem}
\begin{proof}

For $0 \leq r \leq \infty$, let $E_r^{p,q} = \frac{Z_r^{p,q}}{Z_{r-1}^{p+1,q-1} + B_{r-1}^{p,q}}$
endowed with the quotient bornology.  These are the sheets of the spectral sequence being constructed.
Convergence is guaranteed by the boundedness of $F$.
Let $\eta_r^{p,q} : Z_r^{p,q} \to E_r^{p,q}$
be the projection with kernel $Z_{r-1}^{p+1,q-1} + B_{r-1}^{p,q}$.  

The boundary map $d^{p,q}: Z_r^{p,q} \to Z_r^{p+r,q-r+1}$ induces a bounded differential map
$d^{p,q}_r : E_r^{p,q} \to E_r^{p+r,q-r+1}$ yielding the following commutative diagram.

\[\begin{CD}
Z_r^{p,q} @>{d^{p,q}}>> Z_r^{p+r,q-r+1}\\
@V{\eta_r^{p,q}}VV @VV{\eta_r^{p+r,q-r+1}}V\\
E_r^{p,q} @>>{d^{p,q}_r}> E_r^{p+r,q-r+1}
\end{CD}\]

These definitions have the following consequences, used in what follows
\begin{eqnarray*}
\ker d_r^{p,q} & = & \eta_r^{p,q} ( Z_{r+1}^{p,q} )\\
(\eta_r^{p,q})^{-1}( \im d_r^{p-r,q+r-1} ) & = & B_r^{p,q} + Z_{r-1}^{p+1,q-1}\\
Z_{r-1}^{p+1,q-1} \cap Z_{r+1}^{p,q} & = & Z_r^{p+1, q-1}\\
Z_{r+1}^{p,q} \cap (\eta_r^{p,q})^{-1}(\im d_r^{p-r,q+r-1}) & = & B_r^{p,q} + Z_r^{p+1,q-1}
\end{eqnarray*}

This proof will consist of three steps.  The first step consists in verifying that $(E_r, d_r)$
is a bornological spectral sequence.  The next step is to show that it has the appropriate $E_1$-term.
The final step is ensuring that it has the appropriate $E_\infty$-term.
These steps are carried out in the following lemmas.
\end{proof}

\begin{lem}
On the bornological vector spaces $E_r$ associated to $(A,F,d)$,
there is a bornological isomorphism $E^{p,q}_{r+1} \isom H^{p,q}(E_r, d_r)$.
In particular $(E_r, d_r)$ is a bornological spectral sequence.
\end{lem}
\begin{proof}
Let $\gamma : Z_{r+1}^{p,q} \to H^{p,q}( E_r, d_r )$ be the bounded map given by the composition
\[ Z_{r+1}^{p,q} \stackrel{\eta_r^{p,q}}{\to} \ker d_r^{p,q} 
	\stackrel{\pi}{\to} H^{p,q}(E_r^{*,*},d_r) \]
where $\pi$ is the usual projection onto $H^{p,q}( E_r, d_r ) = \frac{\ker d_r^{p,q}}{\im d_r^{p-r,q+r-1}}$.
Since $\ker \gamma = Z_{r+1}^{p,q} \cap (\eta_r^{p,q})^{-1}(\im d_r^{p-r,q+r-1}) = B_r^{p,q} + Z_r^{p+1,q-1}$,
there is an isomorphism of vector spaces
\[ \frac{Z_{r+1}^{p,q}}{B_r^{p,q} + Z_r^{p+1,q-1}} = E_{r+1}^{p,q} \isom H^{p,q}(E_r, d_r) \]
given by $\gamma' : z + (B_r^{p,q} + Z_r^{p+1,q-1}) \mapsto \gamma(z) + (\im d_r^{p-r,q+r-1})$.
We show that $\gamma'$ is the required bornological isomorphism.  

Let $U$ be a bounded subset of $\frac{Z_{r+1}^{p,q}}{B_r^{p,q} + Z_r^{p+1,q-1}} = E_{r+1}^{p,q}$.  
There is a bounded subset $U'$ of $Z_{r+1}^{p,q}$ such that $\eta_{r+1}^{p,q}(U') = U$, so
$\gamma'( U ) = \eta_r^{p,q}(U') + (\im d_r^{p-r,q+r-1})$.  As $\eta_r^{p,q}$ is a bounded map,
$\eta_r^{p,q}(U')$ is a bounded set in $\ker d_r^{p,q}$ and
$\eta_r^{p,q}(U') + (\im d_r^{p-r,q+r-1})$ is bounded in $H^{p,q}(E_r, d_r )`$.  The boundedness of $\gamma'$ is verified.

Let $\phi : \frac{\ker d_r^{p,q}}{\im d_r^{p-r,q+r-1}} \to \frac{Z_{r+1}^{p,q}}{B_r^{p,q} + Z_r^{p+1,q-1}}$
be given by $z + (\im d_r^{p-r,q+r-1}) \mapsto (\eta_r^{p,q})^{-1}(z) \cap Z_{r+1}^{p,q} + (B_r^{p,q} + Z_r^{p+1,q-1})$.
This is the inverse of $\gamma'$.  Let $U$ be a bounded subset of
$\frac{\ker d_r^{p,q}}{\im d_r^{p-r,q+r-1}}$.  There exists a bounded subset $U'$ of $\ker d_r^{p,q}$ such
that $U' + (\im d_r^{p-r,q+r-1})$ contains $U$ in $\frac{\ker d_r^{p,q}}{\im d_r^{p-r,q+r-1}}$.
As $\ker d_r^{p,q} \subset E_r^{p,q}$, $U'$ is bounded in $E_r^{p,q}$,  so there is a bounded subset $U''$
of $Z_r^{p,q}$ with $U' = \eta_r^{p,q}( U'' )$.  Thus $U'' + B_{r-1}^{p,q} + Z_{r-1}^{p+1,q-1}$ is the full preimage of $U'$ under $\eta_r^{p,q}$.
\begin{eqnarray*}
(\eta_r^{p,q})^{-1}(U') \cap Z_{r+1}^{p,q} & = & U''\cap Z_{r+1}^{p,q} + B_{r-1}^{p,q} \cap Z_{r+1}^{p,q} + 
						Z_{r-1}^{p+1,q-1} \cap Z_{r+1}^{p,q}\\
	& = & U''\cap Z_{r+1}^{p,q} + B_{r-1}^{p,q} + Z_r^{p+1,q-1}\\
	& \subset & U'' \cap Z_{r+1}^{p,q} + B_r^{p,q} + Z_r^{p+1,q-1}
\end{eqnarray*}

Thus 
\[\phi(U) \subset U''\cap Z_{r+1}^{p,q} + (B_r^{p,q} + Z_r^{p+1,q-1}) \] 
in $\frac{Z_{r+1}^{p,q}}{B_r^{p,q} + Z_r^{p+1,q-1}}$.
As $U'' \cap Z_{r+1}^{p,q}$ is bounded in $Z_{r+1}^{p,q}$, $\phi(U)$ is bounded in 
$\frac{Z_{r+1}^{p,q}}{B_r^{p,q} + Z_r^{p+1,q-1}}$, whence $\phi$ is a bounded map.
\end{proof}

\begin{lem}
The bornological spectral sequence $(E_r,d_r)$ associated to $(A,F,d)$ has the property that 
$E_1^{p,q} \isom H^{p+q}(F^{p}A/F^{p+1}A)$ as bornological vector spaces.
\end{lem}
\begin{proof}
Since $Z_{-1}^{p+1,q-1} = F^{p+1}A^{p+q}$, 
$B_{-1}^{p,q}=d(F^{p+1}A^{p+q-1})$, and $Z_0^{p,q} = F^pA^{p+q} \cap d^{-1}(F^pA^{p+q+1})$, we have
\begin{eqnarray*}
E_0^{p,q} & = & \frac{Z_0^{p,q}}{Z_{-1}^{p+1,q-1} + B_{-1}^{p,q}} \\
 & = & \frac{F^pA^{p+q} \cap d^{-1}(F^pA^{p+q+1})}{F^{p+1}A^{p+q} + d(F^{p+1}A^{p+q-1})}\\
 & = & \frac{F^pA^{p+q}}{F^{p+1}A^{p+q}}
\end{eqnarray*}

The map $d_0^{p,q} : E_0^{p,q} \to E_0^{p,q+1}$ is induced by $d^{p,q}:F^pA^{p+q} \to F^pA^{p+q+1}$,
fitting into a commutative diagram
\[\begin{CD}
F^pA^{p+q} @>{d}>> F^pA^{p+q+1}\\
@V{\pi}VV @VV{\pi}V\\
E_0^{p,q} = \frac{F^pA^{p+q}}{F^{p+1}A^{p+q}}  @>>{d_0}> \frac{F^pA^{p+q+1}}{F^{p+1}A^{p+q+1}} = E_0^{p,q+1}
\end{CD}\]

where $\pi$ are the usual projections.  
As $H^{p,q}(E_0, d_0)$ is the homology of the complex $(F^pA^* / F^{p+1}A^*, d_0)$,
$H^{p,q}(E_0,d_0) = H^{p+q}(F^pA/F^{p+1}A)$, yielding 
a bornological isomorphism
\[ E_1^{p,q} \isom H^{p+q}(F^pA/F^{p+1}A) \]
\end{proof}

\begin{lem}
Assume the filtration on $(A,F,d)$ is bounded.
The associated bornological spectral sequence $(E_r,d_r)$ converges to 
$H(A,d)$.  That is, \[E_{\infty}^{p,q} \isom F^pH^{p+q}(A,d)/F^{p+1}H^{p+q}(A,d).\]
\end{lem}
\begin{proof}
The filtration $F$ on $A$ induces a
filtration on $H(A,d)$, given by $F^pH(A,d) = \im\{ H(inclusion): H(F^pA) \to H(A) \}$.
Let $\eta_\infty^{p,q} : Z_\infty^{p,q} \to E_\infty^{p,q}$ and $\pi: \ker d \to H(A,d)$
denote the projections.
\begin{eqnarray*}
F^pH^{p+q}(A,d) & = & H^{p+q}( \im( F^pA \to A), d )\\
 & = & \pi( F^pA^{p+q} \cap \ker d )\\
 & = & \pi( Z_\infty^{p,q} )
\end{eqnarray*}

\begin{eqnarray*}
\pi( \ker \eta_\infty^{p,q} ) & = & \pi( Z_\infty^{p+1,q-1} + B_\infty^{p,q} )\\
 & = & \pi( Z_\infty^{p+1,q-1} )\\
 & = & F^{p+1}H^{p+q}(A,d)
\end{eqnarray*}
so $\pi$ induces an isomorphism of vector spaces 
\[ d_\infty : E_\infty^{p,q} \to \frac{F^pH^{p+q}(A,d)}{F^{p+1}H^{p+q}(A,d)}\]

As $\pi: \ker d \to H(A,d)$ is bounded and 
$\pi( Z_\infty^{p,q} ) = F^pH^{p+q}(A,d)$, the restriction  
$\pi : Z_\infty^{p,q} \to F^pH^{p+q}(A,d)$ is a bounded surjection.
Let $U$ be a bounded subset of $E_\infty^{p,q}$.  There is a bounded subset $U'$ of $Z_\infty^{p,q}$
such that $\eta_\infty^{p,q}( U' ) = U$.  As $\pi$ is a bounded map, $\pi(U')$ is a bounded
subset of $F^pH^{p+q}(A,d)$.  Since $d_\infty(U) = \pi(U') + F^{p+1}H^{p+q}(A,d)$ is a bounded subset
of $\frac{F^pH^{p+q}(A,d)}{F^{p+1}H^{p+q}(A,d)}$, $d_\infty$ is a bounded map.

Consider the map
\[\phi : \frac{F^pH^{p+q}(A,d)}{F^{p+1}H^{p+q}(A,d)} \to \frac{Z_\infty^{p,q}}{Z_\infty^{p+1,q-1} + B_\infty^{p,q}}\]
given by $\phi : z + ( F^{p+1}H^{p+q}(A,d) ) \mapsto \pi^{-1}(z)\cap Z_\infty^{p,q} + ( Z_\infty^{p+1,q-1} + B_\infty^{p,q} )$.  This is the inverse of $d_\infty$.
It remains to show that $\phi$ is a bounded map.

Let $U$ be a bounded subset of $\frac{F^pH^{p+q}(A,d)}{F^{p+1}H^{p+q}(A,d)}$.  There is a bounded
$U'$ subset of $F^pH^{p+q}(A,d)$ which projects to $U$.  As
$F^pH^{p+q}(A,d)$ is contained in $H^{p+q}(A,d)$, $U'$ is a bounded subset of $H^{p+q}(A,d)$.
There exists a bounded subset $U''$ in $\ker d^{p+q}$ with $U' = U'' + (\im d^{p+q-1})$.  As $U'$ is a subset of
$F^pH^{p+q}(A,d)$ we can assume $U' \subset \ker d^{p+q} \cap F^pA^{p+q} + (\im d^{p+q-1})$, so
$U'' \subset  Z_\infty^{p,q} + B_\infty^{p,q}$
Therefore $U''$ is bounded in the subspace $Z_\infty^{p,q} + B_\infty^{p,q}$, and
$\pi(U'') = U'' + (\im d^{p+q-1}) \supset U' = U + F^{p+1}H^{p+q}(A,d)$.  Thus 
$\pi^{-1}(U) \subset U'' + \im d^{p+q-1}$.
\begin{eqnarray*}
\pi^{-1}(U) \cap Z_\infty^{p,q} & \subset & U''\cap Z_\infty^{p,q} + \im d^{p+q-1}\cap Z_\infty^{p,q}\\
 & = & U''\cap Z_\infty^{p,q} + B_\infty^{p,q}
\end{eqnarray*}
So
\begin{eqnarray*}
\phi(U) & = & \pi^{-1}(U)\cap Z_\infty^{p,q} + ( Z_\infty^{p+1,q-1} + B_\infty^{p,q} ) \\
 & \subset & U''\cap Z_\infty^{p,q} + ( Z_\infty^{p+1,q-1} + B_\infty^{p,q} )
\end{eqnarray*}
As $U''$ is bounded in $Z_\infty^{p,q} + B_\infty^{p,q}$, $U'' \cap Z_\infty^{p,q}$ is
bounded in $Z_\infty^{p,q}$.  Thus $\phi$ is a bounded map.
\end{proof}

We now move to the application of Lemma \ref{filtDiffGrdBornModSS} in the case which will be of most interest in the sequel.
A double complex of bornological modules is a bigraded module $M = \bigoplus_{p\geq0,q\geq0} M^{p,q}$,
where each $M^{p,q}$ is a bornological module, along with two bounded linear maps $d'$
and $d''$, of bidegree $(1, 0)$ and $(0, 1)$ respectively, satisfying 
$d'^2 =  d''^2 = d'd'' + d''d' = 0$.  The total complex, $(total(M),d)$ of the double
complex $\{M^{*,*}, d', d''\}$ is the differential graded bornological module with 
$total(M)^n = \bigoplus_{p+q=n} M^{p,q}$ and $d = d' + d''$.

There are two standard filtrations on the total complex.
\begin{eqnarray*}
F_I^p( total(M) )^t & = & \bigoplus_{r \geq p} M^{r,t-r}\\
F_{II}^p(total(M))^t & = & \bigoplus_{r \geq p}M^{t-r,r}
\end{eqnarray*}
will be referred to as the columnwise filtration and rowwise filtration respectively.  
Both are decreasing filtrations, respected by the differential.  As $M^{*,*}$ is
first-quadrant, each of these filtrations are bounded and by Lemma \ref{filtDiffGrdBornModSS}
we obtain two spectral sequences of bornological modules converging to $H( total(M),d )$.  
At $M^{p,q}$ there are two boundary maps, $d'$ and $d''$, with respect to each of which
we may calculate a bigraded cohomology of $M$.  Specifically, let 
$H_I^{p,q}(M) = \frac{\im d'': M^{p,q-1} \to M^{p,q}}{\ker d'': M^{p,q} \to M^{p,q+1}}$
and $H_{II}^{p,q}(M) = \frac{\im d': M^{p-1,q} \to M^{p,q}}{\ker d': M^{p,q} \to M^{p+1,q}}$.
In this way, $H^{*,*}_I(M)$ is a double complex with trivial vertical differential and $H^{*,*}_{II}(M)$ is a 
double complexes with trivial horizontal differential.  We may then take cohomology with respect
to the nontrivial boundary map to obtain the iterated cohomology spaces $H^{*,*}_{II} H_I(M)$ and
$H^{*,*}_I H_{II}(M)$ of $M$.

\begin{lem}\label{dCompBornModSS}
Given a double complex $( M^{*,*}, d', d'')$ of bornological modules and bounded maps, there are
two spectral sequences of bornological modules, $( {_I}E_r^{*,*}, {_I}d_r )$ and 
$\{ {_{II}}E_r^{*,*}, {_{II}}d_r \}$ with
${_I}E_2^{p,q} \isom H_I^{*,*} H_{II}(M)$ and ${_{II}}E_2^{p,q} \isom H_{II}^{*,*} H_{I}(M)$.  If 
$M^{*,*}$ is a first-quadrant double complex, both spectral sequences converge to $H^*( total(M), d )$.
\end{lem}
\begin{proof}
The first-quadrant hypothesis is here to ensure convergence of the spectral sequences, and plays
no role in the calculation of the $E_2$-terms.
In the case of $F_I^p$ we have
\[ {_I}E_r^{p,q} = H^{p+q}\left( \frac{F_I^p( total(M) )}{F_I^{p+1}( total(M) )}, d \right) \]
The differential on $total(M)$ is given by
$d = d' + d''$ so that $d'( F_I^p( total(M) ) ) \subset F_I^{p+1}(total(M))$.  There is a
bornological isomorphism
\[ \left( \frac{F_I^p( total(M) )}{F_I^{p+1}( total(M) )} \right)^{p+q} \isom M^{p,q} \]
with the induced differential $d''$, thus ${_I}E_1^{p,q} \isom H_{II}^{p,q}(M)$.

Consider the following maps
\begin{eqnarray*}
  & &i : H^n( F_I^p ) \to H^n(F_I^{p-1})\\
  & &j : H^n( F_I^p ) \to H^n( F_I^p/F_I^{p+1} )\\
  & &k : H^n( F_I^p/F_I^{p+1} ) \to H^{n+1}(F_I^{p+1})\\
  & &d_1 : H_{II}^{p,q}(M) \to H_{II}^{p+1,q}(M)
\end{eqnarray*}
where $i$ is induced by the inclusion $F_I^{p-1} \to F_I^p$, $j$ is induced by the
quotient map $F_I^p \to F_I^p/F_I^{p+1}$, $k$ is the connecting homomorphism, and
$\partial:F_I^p/F_I^{p+1} \to F_I^{p+1}/F_I^{p+2}$ is induced by the differential $d$.
It is clear that $i$ and $j$ are bounded.  The $k$ map sends $[x + F_I^{p+1}] \in H^n(F_I^p/F_I^{p+1})$
to $[dx] \in H^{n+1}(F_I^{p+1})$.  If $U$ is a bounded subset of $H^n(F_I^p/F_I^{p+1})$ then there is
a bounded subset $U'$ in the kernel of $\partial:F_I^p/F_I^{p+1} \to F_I^{p+1}/F_I^{p+2}$ with
$U' + (\im \partial) = U \in H^n(F_I^p/F_I^{p+1})$.  There is $U''$ a bounded subset of $F_I^p$ with
$U' = U'' + F_I^{p+1} \in F_I^p/F_I^{p+1}$.  As $d$ is a bounded map, $d(U'')$ is a bounded subset of $F_I^{p+1}$.
It follows that $[d(U'')]$ is bounded in $H^{n+1}(F_I^{p+1})$, and $k$ is a bounded map.

A class in $H^{p+q}(F_I^p/F_I^{p+1})$ can be written as $[x + F_I^{p+1}]$, where $x \in F_I^p$ and
$dx \in F_I^{p+1}$, or it can be written as a class $[z] \in H_{II}^{p,q}(M)$, $z \in M^{p,q}$.
$k$ sends $[x + F_I^{p+1}]$ to $[dx] \in H^{p+q+1}(F_I^{p+1})$.  Taking $z$ as a representative 
this determines $[d'z] \in H^{p+q+1}(F_I^{p+1})$, since $d''(z) = 0$.  Thus we can
consider $d'z$ as an element of $M^{p+1,q}$.
The map $j$ assigns to a class in $H^{p+q+1}(F_I^{p+1})$ its representative in 
$H^{p+q+1}(F_I^{p+1}/F_I^{p+2})$.   This gives $d_1 = j \circ k$ as the induced mapping of
$d'$ on $H_{II}^{p,q}(M)$, so $d_1 = \bar{d'}$.  Thus 
${_I}E_2^{p,q} \isom H_{I}^{p,q}H_{II}^{*,*}(M)$.
Symmetry gives ${_{II}}E_2^{p,q} \isom H_{II}^{p,q}H_{I}^{*,*}(M)$.
\end{proof}

In the sequel, it will be necessary for us to compare spectral sequences of bornological spaces.
\begin{defn}
Let $(E_r, d_r)$ and $(E'_r, d'_r)$ be two bornological spectral sequences.  A map
of bornological spectral sequences is a family of bigraded bounded linear maps $f = ( f_r : E_r \to E'_r )$,
each of bidegree $(0,0)$, such that for all $r$, $d'_r f_r = f_r d_r$ and $f_{r+1}$ is the
map induced by $f_r$ in cohomology.
\end{defn}

\begin{lem}\label{lem:SSMappingTheorem}
Suppose $f = ( f_r : E_r \to E'_r )$ is a map of bornological spectral sequences, each $E'_r$ is convex and complete, and
each $E_r$ is convex and has a net.  If $f_t$ is a bornological isomorphism for some $t$, then $f_r$ is a bornological
isomorphism for all $r \geq t$.  Moreover, $f$ induces an isomorphism $E_\infty \to E'_\infty$.
\end{lem}
\begin{proof}
It is well known that these isomorphisms exist between the vector spaces.  It remains to show that these
vector space isomorphisms are bornological isomorphisms.  This follows from Lemmas \ref{subspaceNet}, \ref{quotientNet}, 
and Theorem \ref{bddMapIsom}.
\end{proof}

\section{Polynomially bounded cohomology}\label{pcoh}
Let $G$ be a discrete group.
A length function on $G$ is a function $\ell : G \to [0, \infty)$ such that
\begin{enumerate}
\item $\ell(g) = 0$ if and only if $g = 1_G$ is the identity element of $G$.
\item For all $g \in G$, $\ell(g) = \ell(g^{-1})$.
\item For all $g$ and $h \in G$, $\ell(gh) \leq \ell(g) + \ell(h)$.
\end{enumerate}
To a finite generating set $S$ of $G$, we associate a length function $\ell_S$ defined by
$\ell_S( g ) = \min\{ n \, | \, g = s_1 s_2 \ldots s_n \text{ where } s_i \in S \cup S^{-1} \}$.
This length function depends on $S$, but for different choices of $S$ we obtain linearly equivalent
length functions.  We refer to any length function obtained in this way as a word-length function.

Fix some length function $\ell$ on $G$.
For each positive integer $k$ and for $i = 1$ and $i = 2$ define norms on the set of functions $\phi : G \to \C$ by
\[\| \phi \|_{i,k} = \left(\sum_{g \in G} |\phi(g)|^i \left( 1 + \ell(g) \right)^{ik} \right)^{1/i}\]
Let $\mathcal{S}_{\ell}G$ be the set of all functions $f : G \to \C$ such that for all $k$,
$\| f \|_{1,k} < \infty$.  $\mathcal{S}_{\ell}G$ is a Fr\'echet algebra in this family of norms, giving the structure
of a bornological algebra.  In what follows we are solely interested in the case of a word-length
function on $G$.  In this case we denote $\mathcal{S}_{\ell}G$ by $\SG$.
Polynomially equivalent length functions yield the same $\mathcal{S}_{\ell}G$ algebra, so
the particular word-length function used is irrelevant.  If $R$ is a subset of $G$, we also
define $\SR$ to be the subspace of $\SG$ consisting of functions supported on $R$.

Let $A$ be a bornological algebra.  A bornological $A$-module is a complete convex bornological 
space, equipped with a jointly bounded $A$-module structure.  A bornological $A$-module is
bornologically projective if it is a direct summand of bornological module of the form
$A \btensor E$ for some bornological vector space $E$, with the left-action given by the
multiplication in $A$.  An important property of bornologically free modules is that
\[\bHom_{A}(A \btensor B, C) \isom \bHom(B, C).\]

\begin{defn}
The \underline{polynomially bounded cohomology} of $G$ with coefficients in a bornological $\SG$-module $M$ is given by
$HP^*(G; M) = \bExt^*_{\SG}(\C; M)$.
\end{defn}
Here $\bExt$ is the $\Ext$ functor in the bornological category. 
Notice that each of the $\bExt$ groups is a complex bornological vector space.  
Meyer shows in \cite{M1} that this is equivalent to the formulation described in the introduction,
when the coefficient module is $\C$ endowed with the trivial $\SG$-action.
Using $\Ext$ over the topological category one recovers Ji's original definition \cite{Ji1}, 
however from Meyer's work, for trivial coefficients $\C$, the topological and bornological theories coincide. 
  
There is a comparison homomorphism $HP^*(G) \to H^*(G)$ induced by the inclusion $\C[G] \to \SG$. 
An important question with applications to the Novikov conjecture, as well as the $\ell^1$-Bass 
conjecture (see \cite{JOR}) is, "When is this comparison homomorphism is an isomorphism?"  In
\cite{M1} Meyer shows that this is the case for any group equipped with a polynomial 
length combing.  This wide class includes the word-hyperbolic groups of Gromov \cite{G}, 
the semihyperbolic groups of Alonso-Bridson \cite{AB}, and the automatic groups of 
\cite{CEHLPT}, however it does not include all finitely generated groups.  There are examples
of finitely generated groups for which it is known to fail, \cite{JOR2}.

\begin{defn}
A group $G$ is \underline{isocohomological for $M$}, for $M$ an $\SG$-module, if the comparison homomorphism
$HP^*(G;M) \to H^*(G;M)$ is a bornological isomorphism.  It is \underline{strongly isocohomological} if it is isocohomological
for all $\SG$-module coefficients.
\end{defn}
As we will not be interested in weak isocohomologicality, we drop the adjective ``strongly'', and refer to a
group as being isocohomological if it satisfies this strong isocohomologicality condition.

Let $H$ and $Q$ be finitely generated discrete groups with word-length functions
$\ell_H$ and $\ell_Q$ respectively, and let
\[ 0 \to H \stackrel{\iota}{\to} G \stackrel{\pi}{\to} Q \to 0 \]
be an extension of $Q$ by $H$, with word-length function $\ell$.  
( In considering $H$ as a subgroup of $G$, we   
omit the $\iota$ when considering $h \in H$ as an element of $G$. )
Let $q \mapsto \overline q$ be a cross section of $\pi$.
To this cross section there is associated a function, called the factor set of the extension, given by
$[ \cdot, \cdot ] : Q \times Q \to H$ by the formula 
$\overline q_1 \overline q_2 = \overline{q_1 q_2} [ q_1, q_2 ]$.
The factor set has polynomial growth if there exist constants $C$ and $r$ such that 
$\ell_H ([q_1, q_2])  \leq C( (1 + \ell_Q(q_1))( 1 + \ell_Q(q_2)) )^r$.
The cross section also determines a set-theoretic action of $Q$ on $H$ given by
$h^q = \overline{q}^{-1} h \overline{q}$.  
The action is polynomial if
there exist constants $C$ and $r$ such that
$\ell_H( h^q ) \leq C \ell_H(h)  ( 1 + \ell_Q(q) )^r$.

\begin{defn}
An extension $G$ of a finitely generated group $Q$ by a finitely generated group $H$ is said to 
be a \underline{polynomial extension} if there is some cross section yielding a factor set of 
polynomial growth and inducing a polynomial action of $Q$ on $H$.
\end{defn}

An important consequence of this definition is that the word-length function on $H$
is polynomially equivalent to the word-length function on $G$ restricted to $H$.
The following follows from Lemma 1.4 of \cite{N}, and ensures that $\mathcal{S}_{\ell_H}H = \mathcal{S}_{\ell|_H}H$.
\begin{lem}\label{QuasiconvexH}
Let $G$ be a polynomial extension of the finitely generated group $Q$ by the finitely generated group $H$.
There exists constants $C$ and $r$ such that for all $h \in H$,
$\ell(h) \leq \ell_H(h) \leq C(1 + \ell(h) )^r$.
\end{lem}

\begin{lem}
As bornological $\SH$-modules $\SG \isom \SH \btensor \SGH$, where $H$ is endowed with the
restricted length function and $G/H$ is given the minimal length function, 
$\ell^*(gH) = \min_{h \in H} \ell( gh )$, where $\ell$ is the length function on $G$.
\end{lem}
\begin{proof}
Let $R$ be a set of minimal length representatives for right cosets.  Let $r : G \to R$
be the map assigning to $g$, the representative of $Hg$.  Each $g \in G$ has a unique representation
as $g = h_g r(g)$, for $h_g \in H$ and $r(g) \in R$.  
There is an obvious equivalence between $\SGH$ and $\SR$.
Consider the map $\phi : \SG \to \SH \btensor \SR$ given by $\phi(g) = (h_g)\tensor(r(g))$.
This is the desired bornological isomorphism.
\end{proof}

\begin{cor}
A bornologically projective $\SG$-module is a bornologically projective $\SH$-module by
restriction of the $\SG$-action.
\end{cor}

\begin{cor}
Let $M$ be an $\SG$-module. 
Any bornologically projective $\SG$-module resolution of $M$ is a bornologically projective
$\SH$-module resolution of $M$.
\end{cor}

Consider the following:
\[ \ldots \stackrel{\delta}{\to} \SG^{\btensor n} \stackrel{\delta}{\to} 
	\SG^{\btensor n-1} \stackrel{\delta}{\to} \ldots \stackrel{\delta}{\to} 
	\SG \btensor \SG \stackrel{\delta}{\to} \SG \stackrel{\epsilon}{\to} \C \to 0 \]
where $\delta : \SG^{\btensor n} \to \SG^{\btensor n-1}$ is the usual boundary map given
by 
\[ 
\delta( g_1, \ldots, g_n ) = \sum_{i = 1}^{n} \left( -1 \right)^i (x_1, \ldots, \widehat{x_i}, \ldots, x_n)
\]
and extend by linearity, where the tuple $(g_1, \ldots, g_n)$ represents the 
elementary tensor $g_1 \tensor \ldots \tensor g_n$.  As defined, $\delta$ is a bounded map and the map
$s : \SG^{\btensor n} \to \SG^{\btensor n+1}$ given on generators by
\[ s( g_1, \ldots, g_n ) = ( 1_G, g_1, \ldots, g_n ) \] is a
bounded $\C$-linear contracting homotopy for this complex.  This is a bornologically projective
resolution of $\C$ over $\SG$, which we call the \emph{standard bornological resolution} for the group $G$.  

For groups with additional finiteness conditions, there are resolutions with better properties to consider.  By
\cite{JR}, if an isocohomological group $Q$ is of type $HF^\infty$, then there is a bornological
projective resolution of $\C$ over $\SQ$ of the form
\[ \ldots \to R_p \to R_{p-1} \to \ldots \to R_0 \to \C \to 0 \]
with each $R_p$ a bornologically free $\SQ$ module of finite rank.

\begin{thm}\label{MainThm}
Let $0 \to H \to G \rightarrow Q \to 0$ be a polynomial extension of the $HF^\infty$ group $Q$,
with both $H$ and $Q$ isocohomological.
There is a bornological spectral sequence with $E^{p,q}_2 \isom HP^p( Q; HP^q(H) )$ 
which converges to $HP^*(G)$.
\end{thm}

Assuming Theorem \ref{MainThm} and Corollary \ref{isoExtension}, we begin by verifying Corollary \ref{novExtension}.

A group $G$ acts on $\ell^2(G)$ via $(g \cdot f)(x) = f( g^{-1}x)$.  This action extends by linearity
to yield an action by $\C G$ on $\ell^2(G)$ by bounded operators.  The completion of $\C G$
in $\mathcal{B}( \ell^2(G) )$, the space of all bounded operators on $\ell^2(G)$ endowed with
the operator norm, is the reduced group $C^*$-algebra, $C^*_rG$.
Let $\mathcal{S}^2G$ be the set of all functions $f : G \to \C$ such that for all $k$,
$\| f \|_{2,k} < \infty$.  The group $G$ is said to have the Rapid Decay property
if $\mathcal{S}^2G \subset C^*_rG$, \cite{Jol}.  We use the following result of Noskov.
\begin{thm}[\cite{N}]\label{RDExtension}
Let $G$ be a polynomial extension of the finitely generated group $Q$ by the finitely generated group $H$.
If $H$ and $Q$ have the Rapid Decay property, so does $G$.
\end{thm}

\begin{proof}[Proof of Corollary \ref{novExtension}]
By Corollary \ref{isoExtension}, $G$ is isocohomological for $\C$.  By Noskov, $G$ has the
Rapid Decay property.  The result follows from appealing to Connes-Moscovici.
\end{proof}

\section{Proof of Theorem \ref{MainThm}}
Throughout this section, we assume the hypotheses of Theorem \ref{MainThm}.
Let $(P_*, d_P)$ be the standard bornological resolution for $G$, and
let $T_*$ be the tensor product of $P_*$ by $\C$ over
$\SH$.  The polynomial extension properties give $T_q \isom \SQ \btensor \SG^{\btensor q}$.
As the $P_q$ are bornological $\SG$-modules, they are by restriction, bornological $\SH$-modules.
The quotient group $Q$ acts on $\bHom_{\SH}(P_q, \C)$ via 
$\left(q \phi\right)(x) = \overline{q} \cdot \phi\left( {\overline{q}}^{-1} x \right)$,
where $\overline{ \cdot }: Q \to G$ is a cross-section giving the polynomial extension properties.
This extends to a bornological $\SQ$-module structure on $\bHom_{\SH}(P_q, \C)$.
Let $(R_*,d_R)$ be a bornologically projective resolution of $\C$ over $\SQ$ with each $R_p$ finite rank.

Set $C^{p,q} = \bHom_{\SQ}( R_p \btensor T_q, \C ) \isom \bHom_{\SQ}(R_p, \bHom( T_q, \C) )$.  The boundary maps
$d_T$ and $d_R$ induce maps $\delta_T : C^{p,q} \to C^{p,q+1}$ and 
$\delta_R : C^{p,q} \to C^{p+1,q}$ as follows.
\begin{eqnarray*}
\left(\delta_T f\right)(r)(x) & = & (-1)^p f(r)(d_T x)\\
\left(\delta_R f\right)(r)(x) & = & f(d_R r)(x)
\end{eqnarray*}

%

Filter the double complex $C^{*,*}$ by rows.  For a fixed $q$ we have the complex
\[ \ldots \stackrel{\delta_R}{\to} C^{*-1,q} \stackrel{\delta_R}{\to} 
	C^{*,q} \stackrel{\delta_R}{\to} C^{*+1,q} 
	\stackrel{\delta_R}{\to} \ldots \]

The bounded homotopy for the complex $R_*$ induces a contraction on $C^{*,q}$,
so that $E^{p,q}_1 = 0$ for $p \geq 1$ and $E^{0,q}_1 = \bHom_{\SQ}( T_q , \C )$.
The adjointness property gives a bornological isomorphism $\bHom_{\SQ}( T_q,  \C ) \isom \bHom_{\SG}( P_q, \C )$.
This identifies $E^{0,q}_1 \isom \bHom_{\SG}(P_q, \C)$.  As $P_*$ was
a projective $\SG$-complex, we obtain that the $E_2$-term is precisely
$HP^*( G )$, and the spectral sequence collapses here.

%

We now examine the double complex when filtered by columns.
For a fixed $p$ we have the complex
\[ \ldots \stackrel{\delta_T}{\to} C^{p,*-1} \stackrel{\delta_T}{\to} 
	C^{p,*} \stackrel{\delta_T}{\to} C^{p,*+1} 
	\stackrel{\delta_T}{\to} \ldots \] 

By adjointness, $C^{p,q} \isom \bHom_{\SQ}( R_p, \bHom( T_q, \C ) )$, and the boundary map $d_T$
induces a map $d_T^* : \bHom(T_q, \C) \to \bHom( T_{q+1}, \C )$.

\begin{lem}
There are identifications $\ker \delta_T = \bHom_{\SQ}(R_p, \ker d^*_T)$ and 
$\im \delta_T = \bHom_{\SQ}(R_p, \im d^*_T)$.
\end{lem}

\begin{proof}
If $\varphi \in \ker \delta_T$, then for all $r \in R_p$, $(\delta_T \varphi)(r)(x) = 0$ for all 
$x \in T_*$.  
Thus $\varphi(r) \in \ker d^*_T$ for all $r$ and $\ker \delta_T \subset \bHom_{\SQ}(R_p, \ker d^*_T)$.
If $\xi \in \bHom_{\SQ}(R_p, \ker d^*_T)$ then for all $r \in R_p$, $d^*_T \xi(r) = 0$.  That is
$\xi(r)(d_T x) = 0$ for all $x \in T_*$, and $\delta_T \xi$ is the zero map, establishing
$\bHom_{\SQ}(R_p, \ker d^*_T) \subset \ker \delta_T$.

For $\varphi \in \im \delta_T$, there is $f \in \bHom_{\SQ}( R_p, \bHom(T_q, \C) )$ with $\delta_T f = \varphi$.
Thus, for all $r \in R_p$, $\varphi( r ) = (-1)^p d_T^*( f( r ) ) \in \im d_T^*$.  
In particular, $\im \delta_T \subset \bHom_{\SQ}( R_p, \im d_T^* )$.

That $\bHom_{\SQ}(R_p, \im d^*_T) \subset \im \delta_T$ follows from the finiteness condition on $Q$.  Specifically, we use
the bornological isomorphism $\bHom_{SQ}( R_p, M ) \isom \bHom( \overline{R_p}, M )$ for any $\SQ$-module $M$.  Since $\overline{R_p}$ is finite dimensional, any linear map $\overline{R_p} \to M$ is bounded.
For each $\xi \in \im d^*_T$, pick a $\sigma(\xi) \in \bHom(T_q, \C)$ for which $d_T^*( \sigma(\xi) ) = \xi$.  
We do not require $\sigma$ to be a bounded map.  Let $\mathcal{R}$ be a finite basis for $\overline{R_p}$.

Let $\varphi \in \bHom_{\SQ}(R_p, \im d^*_T ) \isom \bHom( \overline{R_p}, \im d^*_T )$.  Define
a map $f : \overline{R_p} \to \bHom( T_q, \C )$ by setting $f( r ) = (-1)^p \sigma( \varphi( r ) )$
for $r \in \mathcal{R}$ and extending by linearity.  This defines a map $f \in \bHom( \overline{R_p}, \bHom(T_q, \C ) )$.
For $r \in \mathcal{R}$, $\delta_T f (r) = (-1)^p d^*_T( f(r) ) = \varphi(r)$.
Thus $\varphi \in \im \delta_T$.
\end{proof}

\begin{lem}\label{lemQuot}
As bornological vector spaces, $\bHom_{\SQ}( R_p, \frac{\ker d^*_T}{\im d^*_T}) \isom \frac{\bHom_{\SQ}( R_p, \ker d^*_T)}{\bHom_{\SQ}( R_p, \im d^*_T)}$ 
\end{lem}
\begin{proof}
Denote by $v$ the map $\frac{\bHom_{\SQ}(R_p, \ker d_T^*)}{\bHom_{\SQ}(R_p, \im d_T^* )} \to \bHom_{\SQ}(R_p, \frac{\ker d_T^*}{\im d_T^*})$,
given by $v( f + \bHom_{\SQ}(R_p, \im d_T^* ))(r) = f(r) + \im d_T^*$.

For $\overline{R_p}$ as above
\begin{eqnarray*}
\bHom(\overline{R_p}, \ker d^*_T) & \subset & \bHom( \overline{R_p}, \bHom( T_q, \C) ) \\
 & \isom & \bHom( \overline{R_p} \btensor T_q, \C ).
\end{eqnarray*}
As $\overline{R_p} \btensor T_q$ is a Fr\'echet space, by Lemma \ref{homNet}, 
Lemma \ref{subspaceNet}, and Lemma \ref{quotientNet},
$\frac{\bHom( \overline{R_p}, \ker d^*_T)}{\bHom( \overline{R_p}, \im d^*_T)}$ has a net.  Moreover,
$\bHom( \overline{R_p}, \ker d^*_T)$ is a complete bornological space.  

The cohomology of the complex 
\begin{eqnarray*}
\bHom( T_q, \C ) & \isom & \bHom( \SQ \btensor \SG^{\btensor q}, \C )\\
	& \isom & \bHom_{\SH}( \SH \btensor  \SQ \btensor \SG^{\btensor q}, \C )\\
	& \isom & \bHom_{\SH }( P_q, \C )
\end{eqnarray*}
is precisely $HP^*(H)$, the polynomially bounded cohomology of the subgroup $H$.
As $H$ is isocohomological for $\C$, by Theorem 11 of \cite{JOR2}, for each $* \geq 0$, $HP^*(H)$ is a finite dimensional 
complex vector space equipped with the fine bornology.  Let $\{ \gamma_1, \ldots, \gamma_k \}$ be a basis for $HP^*(H)$.
For each $\gamma_i$, take an $f_i \in \ker d_T^*$ with $f_i + \im d_T^* = \gamma_i$.  The assignment $\gamma_i \mapsto f_i$ 
extends to a bounded linear map $HP^*(H) \isom \frac{\ker d^*_T}{\im d^*_T} \to \ker d^*_T$ which splits the
quotient map $\ker d^*_T \to \frac{\ker d^*_T}{\im d^*_T} \isom HP^*(H)$.

Let $\phi \in \bHom_{\SQ}( R_p, \frac{\ker d_T^*}{\im d_T^*} )$.  Since $R_p$ has finite rank,  there is a 
$\phi' \in \bHom_{\SQ}( R_p, \ker d_T^* )$ making the following diagram commute.

\begin{equation}
\xymatrix{
	& R_p \ar^{\phi}[d] \ar_{\phi'}[ld] \\
	\ker d^*_T \ar[r] & \frac{\ker d_T^*}{\im d_T^*}
}	 
\end{equation}

This shows the map $\alpha: \bHom_{\SQ}(R_p, \ker d_T^* ) \to \bHom_{\SQ}(R_p, \frac{\ker d_T^*}{\im d_T^*} )$ is surjective.  
As $\ker \alpha = \bHom_{\SQ}( R_p, \im d_T^*)$, $v$ is a bounded linear bijection.
The result follows by Theorem \ref{bddMapIsom}.
\end{proof}

\begin{proof}[Proof of Theorem \ref{MainThm}]
When filtering $C^{*,*}$ by columns, by Lemma \ref{lemQuot} the $E^{p,q}_1 \isom \bHom_{\SQ}( R_p, HP^q(H) )$.
Then $E_2 \isom HP^p(Q; HP^q(H) )$.
As this spectral sequence converges to the same sequence as that obtained when
filtering by rows, we have convergence to $HP^{p+q}(G)$.
\end{proof}

\begin{proof}[Proof of Corollary \ref{isoExtension}]
We compare the polynomial growth spectral sequence with the LHS spectral sequence for the group extension.
The inclusions $\C[H] \to \SH$ and $\C[Q] \to \SQ$ induce a mapping of bornological spectral sequences
$E_r \to E'_r$, where $E_r$ is the spectral sequence resulting from Theorem \ref{MainThm},
and $E'_r$ is the usual spectral sequence associated to the group extension.  
Since $Q$ and $H$ are isocohomological, \[HP^p( Q; HP^q(H) ) \isom H^p(Q; H^q(H) ).\]

The two spectral sequences have bornologically isomorphic $E_2$-terms. 
Both $E_r$ and $E'_r$ are complete convex bornological spaces.  Furthermore, the proof of 
Theorem \ref{MainThm}, when combined with Lemmas \ref{quotientNet} and \ref{subspaceNet}, 
shows that $E_r$ has a net. By Lemma \ref{lem:SSMappingTheorem}
they have bornologically isomorphic limits.
\end{proof}

\end{document}